\documentclass[11pt]{amsart}
\usepackage{amsfonts}
\usepackage{amsmath,amsthm}
\usepackage{latexsym}
\usepackage{amssymb}
\usepackage{enumerate}
\usepackage{comment}

\usepackage{hyperref}
\hypersetup{colorlinks,
	citecolor=blue,
	linkcolor=blue,
	urlcolor=blue}

\DeclareFontFamily{U}{wncy}{}
\DeclareFontShape{U}{wncy}{m}{n}{<->wncyr10}{}
\DeclareSymbolFont{mcy}{U}{wncy}{m}{n}
\DeclareMathSymbol{\Sh}{\mathord}{mcy}{"58}

\usepackage[english]{babel}
\usepackage{color}
\usepackage[latin1]{inputenc}

\numberwithin{equation}{section}

\usepackage{cite}

\theoremstyle{plain}
\newtheorem{theorem}{Theorem}[section]
\newtheorem*{theorem*}{Theorem}

\newtheorem{definition}[theorem]{Definition}

\theoremstyle{remark}

\newtheorem*{lem*}{Lemma}

\newtheorem*{remark*}{Remark}
\newtheorem*{NB*}{NB}

\numberwithin{equation}{section}

\newcommand{\be}{\begin{equation}}
\newcommand{\ee}{\end{equation}}

\newcommand{\vet}{\vec\eta}

\newcommand{\veta}{\vec\eta}

\newcommand{\strela}{\rightharpoonup}

\newcommand{\vxi}{\vec\xi}

\newcommand{\om}{\omega}

\newcommand{\IP}{{\bf P}}

\newcommand{\EE}{{\bf E}}

\newcommand{\PP}{{\mathbf P}}

\newcommand{\cD}{{\mathcal D}}

\newcommand{\cF}{{\mathcal F}}
\newcommand{\cK}{{\mathcal K}}

\newcommand{\cP}{{\mathcal P}}

\newcommand{\cE}{{\mathcal E}}

\newcommand{\SSS}{{\frak S}}

\newcommand{\Z}{{\mathbb Z}}

\newcommand{\N}{{\mathbb N}}
\newcommand{\R}{{\mathbb R}}

\newcommand{\XX}{{\frak X}}
\newcommand{\Fd}{F^\dagger}

\title[Equations with randomness]
{  Long-time behaviour  of solutions  for equations with randomness }

\author{Sergei Kuksin}
\address{ Universit\'e Paris Cit\'e and Sorbonne Universit\'e, CNRS, IMJ-PRG, F-75013 Paris, France 
 \&  Steklov Mathematical Institute of RAS, 
119991 Moscow, Russia \&
 Peoples' Friendship University of Russia (RUDN University),
117198, Russia  \& National Research University Higher School of
	Economics, 119048 Moscow, Russia}
\email{ sergei.kuksin@imj-prg.fr}

\begin{document}
	
\begin{abstract}
We discuss  results on the mixing properties of solutions for discrete and continuous time equations, drived by 
stationary bounded random forces, obtained during the last three years in a number or our joint papers. 
To facilitate applications of these results to specific equations with randomness, the results  are edited and 
 properly re-arranged. 
\end{abstract}

\maketitle

  \section{The equations. }

     Let $(\Omega, \cF, \IP)$ 
     be a probability space, and $H$ and $E$  
      be, accordingly,  a separable Hilbert  and Banach spaces of finite or infinite dimension.
      In what follows,  $H$ is a  phase space and  $E$ is a space of controls for the problems which we examine.
       Let $S:H \times E \to H$ be a continuous map. Consider a discrete time dynamical system in $H$  with added randomness:
 \begin{equation} \label{1}
u_k = S(u_{k-1}, \eta_k^\omega), \quad k\ge1 \qquad u_0 =v; \qquad u_k \in H. 
\ee
Here  $\{\eta_j^\om, j \in \Z\}$ is a random process in $E$. We will denote a trajectory of \eqref{1} as $\{ u_k(v), k\ge0\}$. 
For example,  if $E=H$ and  $S:H\to H$,  then equation \eqref{1} reeds 
 $$
u_k = S(u_{k-1}) + \eta_k^\omega, \quad k\ge1 \qquad u_0 =v.
$$

   Or  let $F: H \times E \to H$ be a mapping with some domain of definition $D(F) \times E$, $D(F) \subset H$, where  
    $D(F)=H$ if dim$H<\infty$,  and $\{\eta^\om(t)\}$ 
   be a random process in $E$, such that the differential equation 
 \begin{equation} 
  \label{2}
\dot u(t) =F(u(t), \eta^\om(t)) \; \quad t\ge0, \qquad u(0) =v
\ee
is well posed, a.s.  We will denote its solution $u(t; v)$. 
 For example,  $E=H$, \ $F:H\supset D(F) \to H$,  and equation \eqref{2} is
 \begin{equation} 
  \label{2'}
\dot u(t) =F(u(t)) + \eta^\om(t) \; \quad t\ge0, \qquad u(0) =v.
\ee
This is a randomly forced differential equation  in $H$. 

 We 
 assume  that process $\eta_k$ in \eqref{1}  is stationary:
for any integer $l$,  probabilistic characteristics of the shifted process $\{ \eta_{k+l}, k\in \Z\}$ are the same as those 
 of $\{ \eta_{k}, k\in \Z\}$.
 Concerning process $\{\eta(t), t\in\R\}$ in \eqref{2}  it will be assumed that it is stationary, or, more generally, 
   statistically $T$-periodic, for some $T>0$. That is,
 probabilistic characteristics of $\eta(t+T)$ are the same as of $\eta(t)$.\footnote{In fact, the stationarity assumption is not necessary for our result which we discuss below; see in \cite{GK26}.}
 
 Our concern is the long-time behaviour of the laws of solutions for \eqref{1} and \eqref{2}, $\cD(u_j)$ and $\cD(u(t))$ (which are Borel measures on $H$).

 \subsection{ Reduction of \eqref{2} to \eqref{1}.} \label{S:1.1}
 It is well known that the problem of long time behaviour of solutions for \eqref{2} may be reduced to that for solutions of some system 
  \eqref{1}. To explain this reduction we firstly suppose that trajectories of  $\eta$ are of low regularity:
 $
 \eta^\om(\cdot) \in L^2_{loc}(\R),
 $
 a.s. As above, the process is assumed to be statistically $T$-periodic. To simplify notation we assume that $T=1$. 
   For $k=1,2,\dots$ we denote $J_k = [k-1, k)$, denote  $J = [0,1)$ and set 
 \be\label{eta}
 \eta_k = \eta\!\mid_{J_k} \in L_2(J_k; E) \sim L_2(J;E) =: \hat E.
 \ee
 The process $\{\eta_k\}$ in space $\hat E$ is stationary.  
 Note that $\hat E$ always   is an infinite-dimensional Hilbert space.

Consider the map
 \[
 S_1: H\times \hat E \to H, \quad (v, \eta(\cdot) ) \to u(1),
 \]
 where $u(t), 0\le t\le1$, solves (2). Then $\{ u_k = u(k)\}$ satisfies 
 \be\label{new_syst}
 u_k = S_1(u_{k-1}, \eta_k), \; \;k\ge1; \;\; \;u_0=v.
 \ee
We thus reduced \eqref{2}  to some system  \eqref{1}.
\medskip

Situation when trajectories of process $\eta$ are regular functions of time is more involved. 
Consider the most important case when 
\be\label{cont}
 \eta^\om(\cdot) \in C(\R),\quad \text{a.s.}
\ee
Now  the discrete time system \eqref{new_syst}
is more complicated than before as  random variables $\eta_k$ and $\eta_{k+1}$ are tied by relation
$
\lim_{t\to1}\eta_k(t) = \eta_{k+1}(0).
$
To handle this difficulty  we will lift system  \eqref{new_syst} to another  one in an extended phase-space. 
Denote ${\cE}= C(\bar J, E)$, $\bar J=[0,1]$,  and 
${\cE_0} =\{ \xi(x)\in {\cE}: \xi(0)=0\}$. Let $\theta(x)$ be the function 
 $\theta(x) = 1-x$, $0\le x\le1$. Then for each $\xi\in {\cE}$, 
$
\xi(x)- \theta(x) \xi(0) \in {\cE_0}.
$
This defines a decomposition of ${\cE}$ to the direct sum 
$$
{\cE} = {\cE_0} \oplus (\theta \otimes Y),
$$
where the corresponding projection $\pi: {\cE} \to {\cE_0}$ sends $\xi(x) \in {\cE}$ to $\xi(x)- \theta(x) \xi(0) \in {\cE_0}$. For  $k\in \Z$ we 
naturally extend $\eta_k$ to a continuous curve on $\bar J_k \sim \bar J$, and  define
$$
\nu_k =\eta(k) = \eta_k(0)=\eta_{k-1}(1)\in E,  \qquad \hat\eta_k = \pi(\eta_k) \in {\cE_0}. 
$$
Then $\eta_k =\nu_k \theta + \hat\eta_k$. We denote $w_k=(u_k, \nu_{k-1}) \in H\times E$, and note that if $\{u_k, k\ge0\}$ is a trajectory of
\eqref{new_syst}, then the sequence $\{w_k, k\ge0\}$, satisfies 
\be\label{ext_syst}
w_k =\hat S_1 (w_{k-1}, \hat\eta_k), \qquad k\ge1, \quad
w_0 =(v, \nu_{-1}).
\ee
Here $\hat S_1$ is the mapping 
$$
\hat S_1:
(H\times E) \times {\cE_0} \to H\times E, \qquad \big((v, \mu), \hat\eta\big) \mapsto 
\big(S_1(v, \mu\theta+ \hat\eta), \hat\eta(1)\big).
$$
The random variables $\{\hat\eta_k\}$  are not tied  any more, and 
 \eqref{ext_syst} turned out to be   an example of a ``good"
system \eqref{1} with phase-space $H:= H\times E$ and space of controls $E:={\cE_0}$.

 \subsection{   Equations with Markovian dynamics.  }
 
 Equation 
 \eqref{1}  defines a Markov process in $H$ (and Markovian dynamics in $H$) if for any $k\ge0$
 \[
\cD(u_l), \; l\ge k+1, \;\; \text{depends only on}\;\; \cD(u_k).
 \]
(In general this is not true, and $\cD(u_l)$ also depends on all $\cD(u_r)$ with $r<l$.) 
In this case,  if $u_0=v\in H$ is non-random, then $\cD(v) = \delta_v$. For $k\ge0$, $\cD(u_k(v))$ depends only
on $v$, and  we denote  
 \[
\cD(u_k(v)) = \cP_k(v; \cdot) \;\text{-- this is  the transition probability from $v$ in time $k$. }
\]
Then for any set $Q\subset H$, \ $\cP_k(v; Q) =   \PP (u_k(v) \in Q)$. Eq.
 \eqref{1} defines a Markovian dynamics if $\{ \eta_k^\om\}$ are i.i.d. random variables. 
 Such  process $\eta$ is called a {\it discrete time white noise in $E$}.

  For \eqref{2} we have similar definitions of   Markovness  and of transition  probabilities 
  $ \cP_t(v; \cdot) $   from $v$ in time $t\ge0$. Eq.
  \eqref{2} defines a Markovian dynamics in $H$ if it has the form \eqref{2'} and there  $\eta^\om(t)$
 is a white noise in $H$,  i.e. $\eta^\om(t) = \dot w^\om(t)$, where $w^\om(t)$ is a Wiener process in space $H$, and the equation is well posed.  More generally, \eqref{2} defines a Markovian dynamics 
 if it has the form
 \be\label{SDE}
\dot u(t) =F(u(t)) +A(u(t)) \dot w(t), \qquad A(u) \; \text{ -- a linear operator in $H$}. 
\ee
Such equations are stochastic ODEs  (SDEs), or  stochastic PDEs (SPDEs)  if dim$\,H=\infty$.
 Vast majority of publications on equations with randomness are about  those with  Markovian 
 dynamics. (But  certainly it is not true that majority of applied problems lead to  such stochastic systems.)

\section{ Long time behaviour  of solutions.   }

   A remarkable feature of eqs. \eqref{1} and \eqref{2} is a highly efficient theory of the limiting behaviour of the laws of their
   solutions  as time goes to infinity.  This theory may be compared with corresponding results for the
     hyperbolic dynamical systems (e.g.  with the results of Sinai and others  on Anosov's diffeomorphisms), but the 
     theory for eqs. \eqref{1},  \eqref{2}   applies to a (different)  much larger class of equations.

   \begin{definition} Equation 
    \eqref{1} is called mixing if there is a measure $\mu$ on $H$ such that for each initial 
    data $v$,
 \begin{equation} 
  \label{mix}
\cD(u_k(v)) \strela \mu \quad\text{as}\quad k\to\infty. 
\ee
   The measure $\mu$ is a {\it stationary measure for \eqref{1}}, and the rate of convergence in \eqref{mix} 
    is called the {\it rate of  mixing}.   
   Definition of mixing for \eqref{2} is the same. 
   \end{definition}
   
   Under some mild additional restrictions on the system (which are satisfied under the assumptions of the theorems which we
   discuss below), for some  random initial date $v$, measurable with respect to the sigma-algebra, generated by the random 
   variables $\{\eta_k, k\le0\}$, the corresponding solution $\{u_k\}$ is such that $\cD (u_k) \equiv\mu$.

 The most common is  the exponential mixing. 
Physicists call the stationary measure $\mu$ {\it a statistical equilibrium} for \eqref{1} (or for \eqref{2}). For a mixing 
system \eqref{1}, if 
  $f$ is a bounded continuous   functional on $H$ (``an observable"),  then
 $\ 
\EE f (u_k(v)) \to \int f(u) \mu(du),
$
for any initial data $v$. -- This explains the importance  of  mixing for physics and  other applications, leading to equations 
\eqref{1} and \eqref{2}.

We call this property of solutions for eqs \eqref{1} or \eqref{2} ``mixing" since solutions $\{u_k\}$ and $u(t)$ 
for such equations are mixing random processes in space $E$, in a  sense very close to the classical:
  if $u(t)$ solves  equation \eqref{2} which is mixing, 
then it holds that  random variables $u(t)$ and $u(t+T)$ are ``almost independent" if $T\gg1$. See  \cite[Section~4.2]{KS25}
and \cite[Section~3.2]{GK26} for discussions.

 \subsection{  Markovian mixing.   }
 
 Mixing in Markovian systems \eqref{1}  and \eqref{2} with dim$\,H <\infty$ is a classical topic, starting the works of Kolmogorov
 and Doeblin, e.g. see in \cite{Bor}. 
 SDEs \eqref{SDE}  are Markovian systems in finite-dimensional phase-spaces. 
  For them  the problem of
mixing  reduces to some questions about linear parabolic PDEs in $\R^n$. Usually for these 
equations stationary measures are of the form $\mu=p(x)dx$, where $p\ge0$ is a non-trivial solution of some 
homogeneous elliptic PDE, see \cite{Kha}. 
 This class of equation with randomness 
 is the only one where stationary measures $\mu$ are relatively well 
understood. 

If in \eqref{2} dim$\,H=\infty$, then we talk about mixing in  SPDEs. A famous example of such equation is 
 the  stochastic  2d\,NSE. 
The mixing in SPDEs has been intensively studied starting the year 2000, see in \cite{KS12}. Still no good general theorem is available for that, 
and even if for some SPDE we may establish the mixing, we know almost nothing about the corresponding 
stationary measure.

 \section{ Non-Markovian mixing.   }
 
As it was explained, equations \eqref{2} may be reduce to systems \eqref{1}. So we will mostly talk  about the latter, 
  where the processes $\eta$ are stationary, but they are not discrete time white noises (i.e.  the random variables 
$\{\eta_k\}$ are not  independent).
Then   solutions $\{u_k\}$ are not Markov processes. A theorem that establishes 
mixing for a large class of equations \eqref{1} was recently proved in \cite{KS25} (also see \cite{GK26}).  
There we assume that the process $\eta$  is stationary and  compactly supported, 
 \begin{equation} 
  \label{bound}
\eta_k\in\cK \qquad \forall\, k, \om,
\ee
where $\cK$ is a compact subset of $E$, 
and that    system \eqref{1} is dissipative in the sense that there exists a compact set
 $K\Subset H$ such that eventually every trajectory $u_k(v)$ enters $K$ and then stays there. We will always suppose  that 
 the map $S$ is $C^2$-smooth. 
There are   more restrictions on  the 
process $\eta$ and mapping $S$ which we will now explain. We will start with systems with finite-dimensional phase 
spaces $H$ 
(but spaces $E$ may be infinite-dimensional).

 \subsection{ Dimension of $H$ is finite.  }
 Consider the ``process of the past":
 \be\label{process}
\vet_k =(\dots, \eta_{k-1}, \eta_k) \in E^{\Z_-}, \qquad k\in\Z,
\ee
where $\Z_-$ is the set of non-positive integers. A classical fact from the theory of probability   is that  
 $\{ \vet_k, k\in \Z\}$ is a stationary Markov process in $E^{\Z_-}$.  By  \eqref{bound}, it also is  a stationary Markov 
process in  $\cK^{\Z_-}$. 
Consider  time-one transition for this process:
 \[
 \cK^{\Z_-}\ni
\vet_k =  (\dots, \eta_{k-1}, \eta_k) \mapsto (\dots, \eta_{k-1}, \eta_k,   \eta_{k+1}) =\vet_{k+1}\in \cK^{\Z_-}.
\]

 Its time-one  transition probability 
 $\cP_1(\vet_k; \cdot)$   is a measure on ${\cK}^{\Z_-}$ of the form 
 \[
\cP_1(\vet_k; \cdot) = \delta_{\vet_k} \otimes Q(\vet_k; \cdot),
\]
where $Q(\vet_k; \cdot)$ is a measure on ${\cK}$.\footnote{The function $Q(\vxi, A)$, which is a Borel measure in 
 $A\subset \cK$, depending on $\vxi\in\cK^{\Z_-}$, is called 
  a {\it transition probability from $\cK^{\Z_-}$ to $\cK$.}}
This  $Q$ is  
the {\it (regular) conditional probability}  for the law of random variable 
$\eta_{k+1}$, given the past $\vet_k$ (see \cite{Dud02}).

In the space $E^{\Z_-}$ we introduce a weighted $\ell_1$-norm:
 \be\label{norm}
 \text{for} \;\; \vet=(\dots, \eta_{-1}, \eta_0), \qquad  \||\vet\|| = {\sum}_{j=0}^\infty \alpha^j |\eta_{-j}| 
\ee
for some fixed $ \alpha \in (0,1)$, 
and furnish  ${\cK}^{\Z_-}$ with a distance, corresponding to this norm.

Let $\cP ({\cK})$ be the set of probability measures on $\cK$,  provided with the total variation distance $\|\cdot\|_{var}$. 

We impose five restrictions on process $\eta$ and map $S$:

   1) {\it Process $\eta$ is exponentially  mixing.} The mapping 
    $$
  {\cK}^{\Z_-} \to \cP({\cK}) , \quad \vet  \mapsto Q(\vet; \cdot)
 $$
 is Lipschitz-continuous. 
 
   We call this property {\it mixing} since it implies that the conditional law of $\eta_{k+1}$ with  fixed 
    past $\vet_k$, \, denoted by us 
 $ Q(\vet_k; \cdot)$,  \ {\it forgets the past exponentially fast}. Indeed, due to the definition of the norm $\||\cdot\||$,
 if  in the vector of  past $\vet_k$
 we change 
 only  its component $\eta_{k-j}$ with some  $j\ge1$,  then the Lipschitz constant of $ Q(\vet_k; \cdot)$ as a function of \ 
 $\eta_{k-j}$ is \ $\le \,$Const$\,\alpha^{-j}$. 
 It can be proved that if a process $\eta$ is mixing in the sense above, then it is mixing in some sense, close to the classical
 (but this not immediate); see \cite[Section~3.2]{GK26}.

   2) {\it Process $\eta$ is recurrent to the origin.}  For any $n\in\N$ and $\delta>0$ there is an integer $s\ge0$ such that 
   $$
   \inf_{\vet \in \cK^{\Z_-} }\IP\{ |\eta_k| <\delta, k=s+1, \dots, s+n \mid \vet_0=\vet\}>0.
   $$

   Let $F$ be a finite-dimensional subspace of $E$ and $\Fd\subset E$ be its complementary subspace, i.e. $E=F\oplus\Fd$ (a space $\Fd$ exists as a consequence of the Hahn-Banach theorem). Accordingly we will write 
   $
   E\ni y=(y_F, y_F^\dagger).
   $
    Let $P_F: E\to F$ and $P_\Fd:E\to \Fd$ be 
   the associated projection operators. Then for any $\vxi\in \cK^{\Z_-}$ the measure $Q(\vxi;\cdot)$ admits a 
   disintegrating with respect to $P_F^\dagger$:
   \be\label{des}
   Q(\vxi;dy) = Q_F^\dagger(\vxi;dy_F^\dagger) Q_F(\vxi, y_F^\dagger;dy_F), \qquad y=(y_F, y_F^\dagger) \in \cK,
   \ee
   where $Q_F^\dagger(\vxi;\cdot)$ is  a projection of  the measure $ Q(\vxi;\cdot)$ to $\Fd$, and 
   $Q_F(\vxi, y_F^\dagger;\cdot)$ is a 
   Borel measure on $\Fd$. E.g. see \cite[Theorem~10.2.1]{Dud02}.

   3)    {\it Process  $\eta$ is regular.}   There exist spaces $F, \Fd$ as above such that in \eqref{des}
   $$
   Q_F(\vxi, y_F^\dagger; dy_F)= \rho_F(\vxi, y_F^\dagger, y_F) dy_F,
   $$
   where $\rho_F$ is a Lipschitz function on $\cK^{\Z_-} \times \cK$, and $dy_F$ is the Lebesgue measure on space $F$. (Space
   $F$ will re-appear below in item 5)).

  4) {\it Free equation \eqref{1}  is globally  dissipative.} If in \eqref{1} $\eta_k \equiv 0$, then 
 $\ 
  u_k(v) \to 0\quad\text{as}\quad k\to\infty, 
   $
 uniformly in  $v\in\cK$. 
 \smallskip
 
 We recall that the mapping $S$ is $C^2$-smooth,  that 
 in \eqref{1}, $\eta_k^\om\in \cK$ \  for all $ k, \om$, 
 and  that by the dissipativity, eventually each trajectory $u_k$ of \eqref{1} 
 stays in $ K \Subset E$.

  5)  {\it Equation \eqref{1} is linearly controllable.} For 
  all $u\in K$, $\eta\in \cK$, the linearised  mapping
 $\ 
 d_\eta S(u, \eta) : F \to H,
 $
 is surjective, where    $F\subset E$ is the space  in 3). 
 \medskip
  
 Note that 
  the  restrictions which we have imposed 
 on  process $\eta$  are separated from those on the mapping $S$, and that the 
 restrictions on    $\eta$ in fact  are those on its transition probability $Q(\vet; \cdot)$.

 So, concerning  process $\eta$ we have assumed that it is 
  (exponentially)   mixing and 
    regular (finite-dimensional components of its conditional distributions have Lipschitz densities). While 
    concerning  map $S$ -- that it is 
   globally dissipative if $\eta=0$, and is linearly controllable \ (the latter is a kind of transversality  condition).

\begin{theorem}\label{t_1}
  Under the imposed assumptions 1)-5) 
   equation \eqref{1}, where $\dim H<\infty$,
    is exponentially mixing. Moreover, for any $v_1, v_2 \in {\cK}$, 
   $
   \| \cD(u_k(v_1)) - \cD(u_k(v_2)) \|_{var} \le C e^{-\gamma k},
   $
   for some positive constants $C$ and $\gamma$.
 \end{theorem}

 {\it Few words about the proof.  }
 Recall that the process $\vet_k$ in $\cK^{\Z_-}$ is Markovian. Let $\{u_k, k\ge0\}$ be a solution of (1). Denote 
 \[
  U_k = (u_k, \vet_k) \in K \times \cK^{\Z_-} =: \XX,
 \]
and consider the map 
 \[
  \SSS : \XX \times \cK \to \XX, \qquad  \big((u, \vet), \xi\big) \mapsto \big( S(u, \xi), (\vet, \xi)\big). 
 \]
Then random trajectory $\{ U_k\}$ satisfies the following relations 
 \begin{equation}\label{33}
U_k = \SSS(U_{k-1}, \eta_k), \quad k\ge1,
\ee
and \ 
 $
 U_0 =(v, \vet_0). 
 $
It turned out that properly understood, system \eqref{33} defines a Markov process in space $\XX$ (this process is rather complicated). 

 So equation  \eqref{33} is a Markovian lifting of  \eqref{1}.

\begin{theorem} \label{t_2}
 The Markov process, defined by  equation \eqref{33} is exponentially mixing. 
\end{theorem}

 This result implies Theorem \ref{t_1}. Indeed, by Theorem \ref{t_2}
   distribution of the pair $U_k = (u_k, \vet_k) $ converges exponentially fast   to some
 measure $M$ on the  space  $ K \times {\cK}^{\Z_-} $. Thus distribution of $u_k$ converges exponentially
  to the projection $\mu$ 
 of measure $M$ to the set  $K$. So \eqref{1} is  exponentially mixing. 
 
 To prove Theorem \ref{t_2} we apply Doeblin's coupling in the form of a  
   method of Kantorovich functional, suggested by the author and his collaborators in 2002 as a  tool to prove the
    mixing for  stochastic PDEs, and  much developed since then; see \cite{KS12} and Introduction in \cite{KS25}.
 The proof of the mixing  uses Newton's quadratic convergence, but in  difference with what we are used to 
 (e.g. in the  KAM theory), 
 now the Newton method
 does not imply super-exponential convergence to a limit --  instead  the power of the method is used to cope with the 
 complexity of the situation. 
 
 \subsection{Dimension of $H$ is finite or infinite.} In this case we still assume  1), 2) and 4) from the previous 
 subsection, and instead of 3) and 5)  suppose that

   $3')$    {\it Process  $\eta$ is regular.}   There is a sequence of subspaces 
   $F_n\subset E$, where $F_1\subset F_2\subset F_3\dots$,
  which either is infinite, or reduces to one space $F_1$. Spaces $F_1, F_2, \dots$ have 
   complementary subspaces $\Fd_1, \Fd_2, \dots$, such that each
  pair $(F_n, \Fd_n)$ satisfies assumption 3). 
  
  We recall that a closed subspace $G\subset H$ is called {\it determining} for system \eqref{1} if there is a constant 
 $ \kappa\in(0,1)$ such that 
 $$
 |(I -P_G)D_uS(u, \eta)|_{H,H} \le \kappa \quad \forall\, u\in K,\ \eta\in \cK,
 $$
 where $P_G$ is the orthogonal projection to $G$.  It is known that 
 system \eqref{1} has a finite-dimensional determining subspace if the map $S$ is smoothing in the sense that for some Hilbert 
 space $V$, compactly and densely embedded in $H$,   $S$ is a $C^1$-smooth map from $H\times E$ to $V$. 
  
   $5')$  {\it Smoothing property and linearised   controllability.} 
System \eqref{1} has a finite-dimensional determining 
   subspace $G$ such that the closure in $H$ of the set $D_\eta S(u, \eta) (F_\infty) $ contains $G$ for all $(u,\eta)$ from 
   some neighbourhood of $K\times \cK$ in $H\times E$. Here $F_\infty =\cup F_j$, where $F_j$ are the spaces from 
   assumption $3')$.

\begin{theorem}\label{t_3}
  Under the assumptions  $1), 2), 3'), 4)$ and $5')$
   equation \eqref{1}, where $\dim H\le \infty$,
    is exponentially mixing. 
     \end{theorem}

 We believe that under the theorem's assumption, in general, equation \eqref{1} is not mixing in the norm of total variation, if
 dimension of $H$ is infinite. 
  Proof of Theorem \ref{t_3} is a generalisation of that of the first assertion of Theorem \ref{t_1}.
  As above, it goes through lifting  the 
 dynamics of equation \eqref{1} to the Markovian dynamics \eqref{33}. 
 
 If in the definition of the norm $\|| \cdot\||$ in 
  \eqref{norm} we replace the exponentially decaying factor $\alpha^{-j}$ with an algebraically decaying one $(1+j)^{-m}$, 
 $m>2$, and if assumption 1) holds with respect to this modified norm, 
  then process $\eta$ is algebraically mixing. In this case, under the hypotheses  of Theorems~\ref{t_1} and \ref{t_3}, 
 system \eqref{1} is algebraically mixing; see \cite{KS27}.

  \section{Where this is proved}
  
  Essentially, Theorem~\ref{t_3} and the first assertion of Theorem~\ref{t_1} were established in \cite{KS25}, while the second assertion of Theorem~\ref{t_1} was derived from the first in \cite{KS26}. In \cite{GK26}, the proofs of both theorems were revised to make them more accessible, and were extended to a broader class of random forces (in particular, 
  the process ${\eta_k}$ was no longer assumed to be stationary). Moreover, the assumptions of the theorems were slightly weakened, allowing them to be applied to a wider class of differential equations \eqref{2}. To simplify the presentation, it was assumed in \cite{GK26} that the space of controls $E$ is Hilbert.

The assumptions imposed on system \eqref{1} in Theorems~\ref{t_1} and~\ref{t_3} are exactly those of \cite{GK26}, except that the space $E$ is now allowed to be Banach rather than Hilbert, as required for some applications below. The proofs in the Banach setting remain essentially unchanged, and any additional technical details can be found, if needed, in \cite{KS25}.

Finally, the reduction of equations \eqref{2} with continuous random processes $\eta(t)$ to ``good'' systems
 \eqref{ext_syst}, presented in Section~\ref{S:1.1}, is an improved version of the construction given in \cite[Section~6.4.2]{GK26}.

  \section{Applications}

Theorems \ref{t_1} and \ref{t_3} naturally apply to discrete time systems \eqref{1} with $\dim H\le\infty$, if the random
process $\eta$ satisfies restrictions $1), 2), 3')$, and $C^2$-smooth mapping $S$ satisfies 4) and $5')$. Conditions 2)
and $3')$ are rather mild restrictions on $\eta$ which demand that the process is regular, and that it 
returns to the vicinity of the 
origin. The most restrictive is condition 1) which requires that the law of $\eta_1$  Lipschitz-continuously
 depends on the past $\vet_0$
 in  the norm of total variation. There are plenty of examples of process $\eta$ which meet these restrictions, see in 
\cite{KS25, GK26}. Assumptions 4) and $5')$ hold for many equations  and are relatively easy to check. 

Among various  equations \eqref{2},  the most important are those, where the process $\eta$ has
continuous trajectories. Let first space $H$ be of finite dimension. Then, if $F$ is $C^2$-smooth and eq.~\eqref{2} is well 
posed, it may be reduced to a $C^2$-smooth system \eqref{ext_syst}. So it has to be checked that the corresponding 
process 
$\hat\eta_k\in {\cE_0}$ and map $\hat S_1$ satisfy assumptions 1)-5). The first three of them, essentially are some 
restrictions on the process $\eta_k$. They are met if process $\eta(t)$ in \eqref{2}
is such that for a segment $J_k$ the law
of the modified restriction $\hat\eta_k$ of $\eta(t)$ to $J_k$ Lipschitz--continuously
 depends on the past in the norm of total 
variation, and some regularity assumption holds. There are plenty of statistically periodic 
bounded processes with continuous trajectories which meet 1)-3), see \cite[Section~6.3]{GK26}, and we believe 
that such  processes are typical among bounded stationary processes with continuous trajectories with the property that 
$\cD(\hat\eta_1)$ Lipschitz--continuously  depends on the past $\vec{\hat\eta}_0$
 in the norm of total variation. 
The map $\hat S_1$ is analysed in \cite{GK26} for the case when process $\eta(t)$ has $H^1$-smooth trajectories, but 
that analyses also applies to the case at hand. It shows that $\hat S_1$  is a well defined $C^2$-smooth mapping,
satisfying assumption 4), if trajectories of $\eqref{2}$ with $\eta\equiv0$ 
converge to the origin, uniformly in $v$ from bounded sets in $H$, and if equation 
$\eqref{2} $ admits a $C^1$-smooth Lyapunov function $V(x)$ which goes to infinity with $|x|_H$ and satisfies 
$$
\big( \nabla V(x), f(x,\eta)\big) \le -\delta \qquad \text{ if $\eta\in\cK$ and $V(x)\ge V_0$,}
$$
for some $V_0, \delta>0$, where $(\cdot, \cdot)$ is the inner product in $H$. If in addition for each $x\in K$ and $\eta\in\cK$ the 
linear operator $D_\eta S(x,\eta): Y \to H$ is surjective, then 5) also holds (the latter restriction holds trivially for
equations \eqref{2'}). 

Concerning Theorem \ref{t_3}, it is shown in \cite{KS25, GK26} that it applies to many equations of mathematical physics 
 with continuous in time statistically-periodic  random forces which 
satisfy  assumptions $1),2),3')$  (there are plenty of  random forces  like that). For example, 
in \cite[Section~6.5]{GK26} we show that the theorem applies to the  randomly forced 3d primitive equations of 
atmospheric mechanics, stirred by a random force as above. 

In Section 3.2 of \cite{GK26} we make the following curious observation. Let us apply Theorem~\ref{t_2} to the  lifted 
Markov system \eqref{33}, corresponding to  trivial system \eqref{1}, where $H$ is the zero-dimensional space 
$\R^0= \{0\}$. Then assumptions $4), 5')$ and $3')$ with $F=\{0\}$ become vacuous, whereas assumptions 1) and 2) 
reduce to conditions on the transition probability $Q(\vxi; \cdot)$. This yields the following corollary to  Theorem~\ref{t_2},
related to Dobrushin's result \cite[Case~2]{Dob} on reconstructing a random process from its conditional 
probabilities: {\it If a transition probability $Q(\vxi; \cdot)$ from $\cK^{\Z_-}$ to $\cK$ satisfies assumptions 1) and 2),
then there exists 	 a stationary process $\eta(t)$, unique in distribution, whose regular conditional distribution 
$
\cD\big(\eta_1\mid \veta_0 
= \vxi\big)
$
equals $Q(\vxi; \cdot)$. Moreover, this process is exponentially mixing in the sense that for any $\vxi\in\cK^{\Z_-}$,  the signed 
measure 
$
\big(\cD(\eta_k \mid  \veta_0 = \vxi) - \cD(\eta_k)\big)
$
weakly converges to zero as $k\to\infty$, exponentially fast. }

\medskip

\noindent
{ \bf Acknowledgements.} The author  acknowledges the support from the Basic Research Program of HSE University.


\begin{thebibliography}{99}
   
\bibitem[Bor98]{Bor}
A.~A. Borovkov, \emph{{Ergodicity and Stability of Stochastic Processes}}, John
  Wiley \& Sons, Ltd., Chichester, 1998.
  
  

\bibitem[Dob70]{Dob}
R.~L. Dobru\v{s}in, \emph{Prescribing a system of random variables by
  conditional distributions}, Teor. Verojatnost. i Primenen. \textbf{15}
  (1970), no.~3, 458--486.

  
\bibitem[Dud02]{Dud02}
R.~M. Dudley, \emph{Real {A}nalysis and {P}robability}, Cambridge University
  Press, Cambridge, 2002.

   \bibitem[GK26]{GK26} P. Gao and S. Kuksin, \emph 
{  Long-time behaviour of  dynamical systems driven by bounded  mixing noises}, arXiv:2607.05981v2.

\bibitem[Kha12]{Kha}
R.~Khasminskii, \emph{{Stochastic Stability of Differential Equations}},  2nd edition,  
  Springer, Heidelberg, 2012.

\bibitem[KS12]{KS12}
S.~Kuksin and A.~Shirikyan, \emph{Mathematics of {T}wo-{D}imensional
  {T}urbulence}, Cambridge University Press, Cambridge, 2012.
   
\bibitem[KS25]{KS25}
S.~Kuksin and A.~Shirikyan, \emph{Mixing for dynamical systems driven by
  stationary noises}, Geom. Funct. Anal
  \textbf{35} (2025),  1346--1399.

\bibitem[KS26]{KS26}
S.~Kuksin and A.~Shirikyan, \emph{Markovian reduction and exponential mixing in total variation for random dynamical systems},
J. Dynam. Differential Equations \textbf{38} (2026), 437-454.


\bibitem[KS27]{KS27}
S.~Kuksin and A.~Shirikyan, \emph{Asymptotics of trajectories for dynamical systems with finitely many equilibria and algebraically mixing noises},
preprint (2026).

   



		
	\end{thebibliography}
 \end{document}